\documentstyle[twoside]{article}
\include{graphicx}
\title{\bf An expansion for the sum of a product of an exponential and a Bessel function}
\author{\sc R. B. Paris\footnote{E-mail address:\ \ {\tt r.paris@abertay.ac.uk}}\\
\\
{\em Division of Computing and Mathematics,}\\
{\em Abertay University, Dundee DD1 1HG, UK}\\
}

\begin{document}
\newcommand{\bee}{\begin{equation}}
\newcommand{\ee}{\end{equation}}
\def\f#1#2{\mbox{${\textstyle \frac{#1}{#2}}$}}
\def\dfrac#1#2{\displaystyle{\frac{#1}{#2}}}
\newcommand{\fr}{\frac{1}{2}}
\newcommand{\fs}{\f{1}{2}}
\newcommand{\g}{\Gamma}
\newcommand{\br}{\biggr}
\newcommand{\bl}{\biggl}
\newcommand{\ra}{\rightarrow}
\renewcommand{\topfraction}{0.9}
\renewcommand{\bottomfraction}{0.9}
\renewcommand{\textfraction}{0.05}
\newcommand{\mcol}{\multicolumn}
\date{}
\maketitle
\pagestyle{myheadings}
\markboth{\hfill {\it R.B. Paris} \hfill}
{\hfill {\it A modified Bessel function sum } \hfill}
\begin{abstract} 
We examine convergent representations for the sum of a decaying exponential and a Bessel function in the form
\[\sum_{n=1}^\infty \frac{e^{-an}}{(\fs bn)^\nu}\,J_\nu(bn),\]
where $J_\nu(x)$ is the Bessel function of the first kind of order $\nu>-\fs$ and $a$, $b$ are positive parameters. By means of a double Mellin-Barnes integral representation we obtain a convergent asymptotic expansion that enables the evaluation of this sum in the limit $a\to 0$ with $b<2\pi$ fixed. A similar result is derived for the sum when the Bessel function is replaced by the modified Bessel function $K_\nu(x)$. The alternating versions of these sums are also mentioned.
\vspace{0.4cm}

\vspace{0.3cm}

\end{abstract}

\vspace{0.5cm}

{\bf 1.  Introduction}
\setcounter{section}{1}
\setcounter{equation}{0}
\renewcommand{\theequation}{\arabic{section}.\arabic{equation}}
We consider the sums
\begin{equation}\label{e11}
S_{J,\nu}(a,b)=\sum_{n=1}^\infty \frac{e^{-an}}{(\fs bn)^\nu}\,J_\nu(bn),\quad\mbox{and}\quad S_{K,\nu}(a,b)=\sum_{n=1}^\infty \frac{e^{-an}}{(\fs bn)^\nu}\,K_\nu(bn),
\end{equation}
where $J_\nu(x)$ is the Bessel function of the first kind and $K_\nu(x)$ the modified Bessel function.
The parameters $a$, $b$ are assumed to be positive and the order $\nu>-\fs$ for $S_{J,\nu}(a,b)$ and $\nu\geq 0$ for $S_{K,\nu}(a,b)$. Convergent expansions are derived for these sums that involve the polylogarithm function $\mbox{Li}_s(x)$ of negative integer order $s$ when $a>b$ and the Riemann zeta function when $a\leq b$; the latter case requires $a, b\in (0,2\pi)$ to secure convergence. The evaluation of the alternating versions of (\ref{e11}) follows in a straightforward manner in terms of these sums with parameters $a$, $b$ and $2a$, $2b$.

Our main interest in the sums in (\ref{e11}) is the limit as $a\to0$ when convergence of the first sum becomes slow.
When $a=0$, the sum $S_{J,\nu}(0,b)$ has been considered by Tri\u ckovi\'c {\it et al.} in \cite{TVS}, where approaches using Poisson's summation formula and Bessel's integral were employed to derive convergent expansions.  
The sums $S_{J,\nu}(0,b)$ and $S_{K,\nu}(0,b)$ have been discussed in \cite{P17} using a Mellin transform approach. With the presence of the exponential factor in the above sums it is found that a double Mellin-Barnes integral is required
to derive an expansion as $a\to0$. The treatment of such integrals has been discussed in \cite[Chapter 7]{PK} in the context of multi-dimensional Laplace-type integrals.

\vspace{0.6cm}

\begin{center}
{\bf 2. \ A series representation for $S_{J,\nu}(a,b)$ valid when $a\geq b$}
\end{center}
\setcounter{section}{2}
\setcounter{equation}{0}
\renewcommand{\theequation}{\arabic{section}.\arabic{equation}}
We consider the sum
\bee\label{e21}
S_{J,\nu}(a,b)=\sum_{n=1}^\infty \frac{e^{-an}}{ (\fs bn)^{\nu}}\, J_\nu(bn)\qquad (a>0, \ b>0),
\ee
where $J_\nu(x)$ denotes the usual Bessel function of real order $\nu>-\fs$. In the case $a=0$, the sum has the value \cite[Theorem 1, (2.8)]{P17}
\bee\label{e22}
S_{J,\nu}(0,b)=\frac{\sqrt{\pi}}{b \g(\nu+\fs)}-\frac{1}{2\g(1+\nu)}\qquad (\nu>-\fs,\ 0<b<2\pi).
\ee

Expanding the Bessel function as a series, we have
\[S_{J,\nu}(a,b)=\sum_{n=1}^\infty e^{-an}\sum_{k=0}^\infty \frac{(-)^k(\fs bn)^{2k}}{k! \g(1+\nu+k)}
=\sum_{k=0}^\infty \frac{(-)^k (\fs b)^{2k}}{k! \g(1+\nu+k)} \sum_{n=1}^\infty e^{-an} n^{2k}\]
upon interchanging the order of summation.
If we introduce the  polylogarithm function $\mbox{Li}_s(x)$ defined by 
$\mbox{Li}_s(z)=\sum_{n=1}^\infty n^{-s}z^n$ \cite[(25.12.10)]{DLMF}, we then obtain
\bee\label{e24}
S_{J,\nu}(a,b)=\sum_{k=0}^\infty \frac{(-)^k (\fs b)^{2k}}{k! \g(1+\nu+k)}\,\mbox{Li}_{-2k}(e^{-a}).
\ee

The polylogarithm function of negative integer order $-n$, $n=0,1, 2, \ldots$ is
\begin{eqnarray*}
\mbox{Li}_{-n}(e^{-a})&=&n! a^{-n-1}\bl\{1+\sum_{k\geq 1}\bl\{\bl(\frac{a}{a+2\pi ki}\br)^{\!n+1}+\bl(\frac{a}{a-2\pi ki}\br)^{\!n+1}\br\}\\
&=&n! a^{-n-1}\bl\{1+2\sum_{k\geq 1}\bl(\frac{a}{\sqrt{a^2+4\pi^2k^2}}\br)^{\!n+1}\!\!\cos [(n+1) \phi_k]\br\},
\end{eqnarray*}
with $\phi_k:=\arctan (2\pi k/a)$.
For large $n$ with $a$ fixed we have
\[\mbox{Li}_{-n}(e^{-a})\sim n! a^{-n-1}\qquad (n\to\infty).\]
The late terms in the sum (\ref{e24}) consequently possess the behaviour 
\[\frac{1}{a}\,\frac{(-)^k \g(2k+1)}{k! \g(1+\nu+k)}\,\bl(\frac{b}{2a}\br)^{\!2k}=\frac{(-)^k \g(k+\fs)}{a\sqrt{\pi}\, \g(1+\nu+k)}\,\bl(\frac{b}{a}\br)^{\!2k}=(-)^kO\bl( k^{-\nu-1/2}\bl(\frac{b}{a}\br)^{\!2k}\br)\]
as $k\to\infty$,
so that the sum (\ref{e24}) converges when $a>b$, and when $a=b$ provided $\nu>-\fs$. Calculations with {\it Mathematica} confirm the validity of this statement, where it should be noted that $\mbox{Li}_s(z)$ is a built-in library function.

\vspace{0.6cm}

\begin{center}
{\bf 3. \ A series representation of $S_{J,\nu}(a,b)$ for $a\to 0$}
\end{center}
\setcounter{section}{3}
\setcounter{equation}{0}
\renewcommand{\theequation}{\arabic{section}.\arabic{equation}}
To deal with the case $a<b$, and in particular to obtain an expansion for $S_{J,\nu}(a,b)$ as $a\to 0$, we employ the Cahen-Mellin integral \cite[p.~89]{PK}
\bee\label{e23a}
e^{-x}=\frac{1}{2\pi i}\int_{c-\infty i}^{c+\infty i} \g(s) x^{-s} ds \qquad(|\arg\,x|<\fs\pi)
\ee
and the Mellin-Barnes integral \cite[(10.9.22)]{DLMF}
\bee\label{e23}
(\fs x)^{-\nu} J_\nu(x)=\frac{1}{2\pi i}\int_{c-\infty i}^{c+\infty i} \frac{\g(s)}{\g(1+\nu-s)}\,(\fs x)^{-2s}ds\qquad (x>0,\ \nu>0),
\ee
where $c>0$ in both cases so that the integration path lies to the right of the poles of the integrand situated at $s=0, -1, -2, \ldots\ $.
Combination of these two results leads to the double Mellin-Barnes integral representation
for $S_{J,\nu}(a,b)$ given by, provided\footnote{If the integation path in (\ref{e23}) is bent back into a loop that encloses the poles of $\g(s)$ with endpoints at infinity in $\Re (s)<0$, then the integral holds without restriction on $\nu$; see \cite[p.~115]{PK}.} $\nu>0$,
\begin{eqnarray}
S_{J,\nu}(a,b)&=&\bl(\frac{1}{2\pi i}\int_{-\infty i}^{\infty i}\br)^{\!2} \frac{\g(s) \g(t)}{\g(1+\nu+s)}\,a^{-t} (\fs b)^{-2s} \sum_{n\geq 1} n^{-2s-t}\, ds\,dt\nonumber\\
&=&\bl(\frac{1}{2\pi i}\int_{-\infty i}^{\infty i}\br)^{\!2} \frac{\g(s) \g(t)}{\g(1+\nu+s)}\,a^{-t} (\fs b)^{-2s} \zeta(2s+t)\, ds\,dt,\label{e31}
\end{eqnarray}
where $\zeta(s)$ is the Riemann zeta function. The integration contours are indented to pass to the right of the pole of the zeta function and those of $\g(s)$ and $\g(t)$.

Let us first consider the poles at $s=-n$, $t=-m$ ($m, n=0,1,2, \ldots$) that result from the gamma functions in the numerator of (\ref{e31}). Displacement of the integration paths to the left produces the formal series
\begin{eqnarray*}
T_1&=&\sum_{m,n\geq0} \frac{(-)^{m+n}}{m! n!}\,\frac{a^m (\fs b)^{2n}}{\g(1+\nu+n)}\,\zeta(-2n-m)\\
&=&\frac{\zeta(0)}{\g(1+\nu)}-\sum_{m,n\geq0}\frac{(-)^n}{(2m+1)! n!}\,\frac{a^{2m+1}(\fs b)^{2n}}{\g(1+\nu+n)}\,\zeta(-2m-2n-1),
\end{eqnarray*}
where even values of $m$ (apart from $m=0$ when $n=0$) do not contribute to the double sum on account of the trivial zeros of $\zeta(s)$ at $s=-2, -4, \ldots\ $. Using the functional relation for $\zeta(s)$ \cite[(25.4.2)]{DLMF}
\[\zeta(s)=2^s\pi^{s-1} \zeta(1-s) \g(1-s) \sin \fs\pi s,\]
the value $\zeta(0)=-\fs$ and the duplication formula $\g(2z)=2^{2z-1} \pi^{-1/2} \g(z)\g(z+\fs)$, we find that
\[T_1=-\frac{1}{2\g(1\!+\!\nu)}+\frac{a}{2\pi^2}\sum_{m,n\geq0}\frac{(-)^n\cos \pi(m\!+\!n)}{(2m\!+\!1)! n!}\,\frac{\g(2m\!+\!2n\!+\!2)}{\g(n\!+\!1\!+\!\nu)}\,\zeta(2m\!+\!2n\!+\!2) \,\bl(\frac{a}{2\pi}\br)^{\!2m}\bl(\frac{b}{2\pi}\br)^{\!2n}\]
\bee\label{e32}
=-\frac{1}{2\g(1\!+\!\nu)}+\frac{1}{\pi}\sum_{m,n\geq0} (-)^m A_{m,n}\,\bl(\frac{a}{2\pi}\br)^{\!2m+1}\bl(\frac{b}{2\pi}\br)^{\!2n},
\ee
where
\bee\label{e32a}
A_{m,n}:=\frac{\g(m+n+1)\g(m+n+\f{3}{2})}{m! n! \g(m+\f{3}{2}) \g(n+1+\nu)}\,\zeta(2m\!+\!2n\!+\!2).
\ee
Application of Stirling's formula for the gamma function, together with the fact that $\zeta(2m+2n+2)=O(1)$ 
for large $m$, $n$, shows that $A_{m,n}$ has only algebraic growth in $m$ and $n$ as $m, n\to\infty$. Consequently the double sum in (\ref{e32}) converges when $a<2\pi$ and $b<2\pi$.

We now deal with the pole of the zeta function when $2s+t=1$ with unit residue. Setting $t=1-2s$, we find the contribution
\[\frac{a^{-1}}{2\pi i}\int_{-\infty i}^{\infty i} \frac{\g(s)\g(1-2s)}{\g(1+\nu+s)} \,\bl(\frac{b}{2a}\br)^{\!2s}ds\]
where the integration path is indented to separate the poles of $\g(s)$ and $\g(1-2s)$.
Evaluation of the residues of $\g(1-2s)$ at $s=\fs+\fs k$, $k=0, 1, 2, \ldots$ leads to
the contribution
\bee\label{e33}
T_2=\frac{1}{b}\sum_{k=0}^\infty \frac{(-)^k \g(\fs k+\fs)}{k! \g(\fs+\nu-\fs k)}\,\bl(\frac{2a}{b}\br)^{\!k}.
\ee
The terms in this sum behave like
\[\frac{(-)^k \g(\fs k+\fs-\nu)}{\sqrt{\pi}\,\g(\fs k+1)}\,\bl(\frac{a}{b}\br)^{\!k} \cos \pi(\nu-\fs k)=
(-)^kO\bl( k^{-\nu-1/2} \bl(\frac{a}{b}\br)^{\!k}\br)\qquad (k\to\infty),\]
so that (\ref{e33}) converges for $a<b$, and for $a=b$ provided $\nu>-\fs$.

By separating terms with even and odd $k$, we can express $T_2$ in terms of well-known special functions to yield
\bee\label{e34}
T_2=\frac{\sqrt{\pi}}{b\g(\nu+\fs)}\bl(1+\frac{a^2}{b^2}\br)^{\!\nu-\frac{1}{2}}\!-\frac{2a}{b^2 \g(\nu)}\,{}_2F_1\bl(1,1-\nu;\frac{3}{2}; -\frac{a^2}{b^2}\br).
\ee
Alternatively, the hypergeometric function can be written in a different form using Euler's transformation \cite[(15.8.1)]{DLMF} to yield
\[T_2=\frac{\sqrt{\pi}}{b\g(\nu+\fs)}\bl(1+\frac{a^2}{b^2}\br)^{\!\nu-\frac{1}{2}}\!-\frac{2a}{(a^2+b^2) \g(\nu)}\,{}_2F_1\bl(1,\fs+\nu;\frac{3}{2};\frac{a^2}{a^2+b^2}\br).\]

To summarise, we have the following expansion.
\newtheorem{theorem}{Theorem}
\begin{theorem}$\!\!\!.$ \ Let $\nu>-\fs$ and $a\leq b$. Then the following convergent expansion holds
\[S_{J,\nu}(a,b)=\frac{\sqrt{\pi}}{b\g(\nu+\fs)}\bl(1+\frac{a^2}{b^2}\br)^{\!\nu-\frac{1}{2}}\!-
\frac{1}{2\g(1+\nu)}-\frac{2a}{b^2 \g(\nu)}\,{}_2F_1\bl(1,1-\nu;\frac{3}{2}; -\frac{a^2}{b^2}\br)\]
\bee\label{e35}
+\frac{a}{2\pi^2}\sum_{m,n\geq0} (-)^m A_{m,n}\,\bl(\frac{a}{2\pi}\br)^{\!2m}\bl(\frac{b}{2\pi}\br)^{\!2n},
\ee
provided $a<2\pi$ and $b<2\pi$. The coefficients $A_{m,n}$ are defined in (\ref{e32a}), .
\end{theorem}
When $a=0$ it is seen that (\ref{e35}) correctly reduces to (\ref{e22}).

The above expansion can be arranged in a form suitable for asymptotic calculations in the limit $a\to0$.
We note that the series representation of the hypergeometric function in (\ref{e34}) is already in the form of a convergent asymptotic expansion as $a\to0$. We find after routine manipulation that
\[T_1=-\frac{1}{2\g(1+\nu)}+\frac{1}{\pi}\sum_{k=0}^\infty (-)^k C_k \bl(\frac{a}{2\pi}\br)^{\!2k+1},\]
where
\bee\label{e36a}
C_k=\sum_{n=0}^\infty A_{k,n}\bl(\frac{b}{2\pi}\br)^{\!2n}\qquad (b<2\pi).
\ee
Consequently we obtain
\begin{theorem}$\!\!\!.$\ Let $\nu>-\fs$ and $a\leq b<2\pi$. Then the following (convergent) asymptotic expansion holds
for $a\to 0$
\[S_{J,\nu}(a,b)=\frac{\sqrt{\pi}}{b\g(\nu+\fs)}\bl(1+\frac{a^2}{b^2}\br)^{\!\nu-\frac{1}{2}}\!-
\frac{1}{2\g(1+\nu)}-\frac{2a}{b^2 \g(\nu)}\,{}_2F_1\bl(1,1-\nu;\frac{3}{2}; -\frac{a^2}{b^2}\br)\]
\bee\label{e36}
+\frac{1}{\pi}\sum_{k=0}^\infty (-)^k C_k\,\bl(\frac{a}{2\pi}\br)^{\!2k+1},
\ee
where the coefficients $C_k$ are defined in (\ref{e36a}).
\end{theorem}

Finally, the alternating version of (\ref{e21}) is given by
\[{\hat S}_{J,\nu}(a,b)=\sum_{n=1}^\infty \frac{(-)^{n-1} e^{-an}} {(\fs bn)^{\nu}}\, J_\nu(bn).\]
It is easily seen that
\[{\hat S}_{J,\nu}(a,b)=S_{J,\nu}(a,b)-S_{J,\nu}(2a,2b),\]
from which the expansion as $a\to 0$ can be obtained from Theorem 2.

\vspace{0.6cm}

\begin{center}
{\bf 4. \ The expansion involving the modified Bessel function}
\end{center}
\setcounter{section}{4}
\setcounter{equation}{0}
\renewcommand{\theequation}{\arabic{section}.\arabic{equation}}
The same procedure can be brought to bear on the sum involving the modified Bessel function $K_\nu(z)$ given by
\bee\label{e41}
S_{K,\nu}(a,b)=\sum_{n=1}^\infty  \frac{e^{-an}} {(\fs bn)^{\nu}}\, K_\nu(bn)\qquad (a>0,\ b>0),
\ee
where we consider only $\nu\geq 0$, since $K_{-\nu}(x)=K_\nu(x)$. We remove from consideration the case when $\nu$ has half-integer values when the above sum reduces to a sum of exponentials; see the appendix. 

From the standard definition $K_\nu(x)=\pi (I_{-\nu}(x)-I_\nu(x))/(2\sin \pi\nu)$ valid for non-integer $\nu$ and the series expansions of $I_{\pm\nu}(x)$, we obtain
\bee\label{e42}
S_{K,\nu}(a,b)=\frac{\pi}{2\sin \pi\nu}\bl\{\sum_{k=0}^\infty \frac{(\fs b)^{2k-2\nu}}{k! \g(1\!-\!\nu\!+\!k)}\,\mbox{Li}_{-2k+2\nu}(e^{-a})-\sum_{k=0}^\infty \frac{(\fs b)^{2k}}{k! \g(1\!+\!\nu\!+\!k)}\,\mbox{Li}_{-2k}(e^{-a})\br\}
\ee
when $a>b$.

Use of the integral representation\footnote{There is an error in the sector of validity in \cite[(10.32.13)]{DLMF} and also in \cite[p.~114]{PK}.} \cite[(10.32.13)]{DLMF}
\[(\fs x)^{-\nu} K_\nu(x)=\frac{1}{4\pi i}\int_{c-\infty i}^{c+\infty i}\g(s) \g(s-\nu) (\fs x)^{-2s} ds\qquad (|\arg\,x|<\fs\pi),\]
where $c>\max\{0,\nu\}$, together with (\ref{e23a}), shows that $S_{K,\nu}(a,b)$ can be expressed as the double Mellin-Barnes integral
\bee\label{e43}
S_{K,\nu}(a,b)=\frac{1}{2}\bl(\frac{1}{2\pi i}\int_{-\infty i}^{\infty i}\br)^{\!2}\g(s)\g(s-\nu)\g(t) a^{-t}(\fs b)^{-2s} \zeta(2s+t)\,ds\,dt,
\ee
where, as in (\ref{e31}), the integration contours are indented to pass to the right of the pole of $\zeta(2s+t)$ and those of the gamma functions. Evaluation of the residues at $s=-n$, $t=-m$ and $s=-n+\nu$, $t=-m$ ($m, n=0, 1, 2, \ldots$) yields, provided $\nu\neq 0, \fs, 1, \f{3}{2}, 2, \ldots\,$,
\begin{eqnarray}
T_1&=&\frac{\pi}{2\sin \pi\nu}\sum_{m,n\geq0}\frac{(-)^ma^m(\fs b)^{2n}}{m! n!}
\bl\{\frac{(\fs b)^{-2\nu}\zeta(-2n\!-\!m\!+\!2\nu)}{\g(n\!+\!1\!-\!\nu)}-\frac{\zeta(-2n\!-\!m)}{\g(n\!+\!1\!+\!\nu)}\br\}\label{e43a}\\
&=&\frac{\pi}{2\sin \pi\nu}\bl\{\frac{1}{2\g(1+\nu)}-\frac{1}{\pi}\sum_{m,n\geq0}(-)^{m+n}\bl\{A_{m,n}\bl(\frac{a}{2\pi}\br)^{\!2m+1}\bl(\frac{b}{2\pi}\br)^{\!2n}\nonumber\\
&&\hspace{8cm}+B_{m,n}\bl(\frac{a}{2\pi}\br)^{\!m}\bl(\frac{b}{2\pi}\br)^{\!2n-2\nu} \br\}\br\}\nonumber
\end{eqnarray}
for $a<2\pi$ and $b<2\pi$. The coefficients $A_{m,n}$ are defined in (\ref{e32a}) and 
\bee\label{e44}
B_{m,n}=\frac{\g(\fs m\!+\!n\!+\!\fs\!-\!\nu) \g(\fs m\!+\!n\!+\!1\!-\!\nu)}{\g(\fs m\!+\!\fs)\g(\fs m\!+\!1) n! \g(n\!+\!1\!-\!\nu)}\,\frac{\zeta(m\!+\!2n\!+\!1\!-\!2\nu)}{\csc \pi(\fs m\!-\!\nu)}.
\ee

To deal with the pole of the zeta function when $2s+t=1$, we obtain from (\ref{e43}) 
\[\frac{1}{4\pi ia} \int_{-\infty i}^{\infty i} \g(s)\g(s-\nu)\g(1-2s) \bl(\frac{2a}{b}\br)^{\!2s}ds,\]
where the integration path is indented to separate the poles of $\g(s)$ and $\g(s-\nu)$ on the left from those of $\g(1-2s)$ on the right ($\nu\neq \fs, 1, \f{3}{2}, 2, \ldots\ $). Evaluation of the residues on the right yields the contribution
\begin{eqnarray}
T_2&=&\frac{1}{2a}\sum_{k=0}^\infty \frac{(-)^k}{k!} \g(\fs k+\fs)\g(\fs k+\fs-\nu) \bl(\frac{2a}{b}\br)^{\!k}\label{e45}\\
&=&\frac{\sqrt{\pi}}{2b} \g(\fs-\nu) \bl(1-\frac{a^2}{b^2}\br)^{\!\nu-\frac{1}{2}}-\frac{a}{b^2} \g(1-\nu)<{}_2F_1\bl(1, 1-\nu;\frac{3}{2};\frac{a^2}{b^2}\br).\nonumber
\end{eqnarray}
The late terms in this expansion possess the behaviour $(-)^k (\fs k)^{-\nu-1/2} (a/b)^k$, so that the series in (\ref{e45}) converges for $a\leq b$ ($\nu\geq 0$).

Then we have the following expansion.
\begin{theorem}$\!\!\!.$\ Let $\nu\neq 0, \fs, 1, \f{3}{2}, 2, \ldots $ and $a\leq b$. Then the following convergent expansion holds 
\[S_{K,\nu}(a,b)=\frac{\sqrt{\pi}}{2b} \g(\fs-\nu) \bl(1-\frac{a^2}{b^2}\br)^{\!\nu-\frac{1}{2}}+\frac{\pi \csc \pi\nu}{4 \g(1\!+\!\nu)}-\frac{a}{b^2} \g(1-\nu)\,{}_2F_1\bl(1, 1-\nu;\frac{3}{2};\frac{a^2}{b^2}\br)\]
\bee\label{e46}
\hspace{2cm}-\frac{1}{2\sin \pi\nu}\sum_{m,n\geq0}(-)^{m+n}\bl\{A_{m,n}\bl(\frac{a}{2\pi}\br)^{\!2m+1}\bl(\frac{b}{2\pi}\br)^{\!2n}
+B_{m,n}\bl(\frac{a}{2\pi}\br)^{\!m}\bl(\frac{b}{2\pi}\br)^{\!2n-2\nu} \br\}
\ee
provided $a<2\pi$ and $b<2\pi$. 
The coefficients $A_{m,n}$ and $B_{m,n}$ are defined in (\ref{e32a}) and (\ref{e44}). 

The double sum appearing in (\ref{e46}) can be written in the alternative (asymptotic) form
\bee\label{e46a}
-\frac{1}{2\sin \pi\nu} \sum_{k=0}^\infty (-)^k D_k \bl(\frac{a}{2\pi}\br)^{\!k},
\ee
where the coefficients $D_k$ are defined by
\[D_{2k}=C'_{2k},\qquad D_{2k+1}=C'_{2k+1}+(-)^{k-1}C_k \qquad (k=0, 1, 2, \ldots)\]
with $C_k$ defined in (\ref{e36a}) and 
\[C'_k=\sum_{n=0}^\infty (-)^n B_{k,n} \bl(\frac{b}{2\pi}\br)^{\!2n-2\nu}\qquad (b<2\pi).\]
\end{theorem}

In the case of integer values of $\nu$ we need to take the limiting value of the double sum $T_1$ in (\ref{e43a}). For example, if $\nu=0$ we have upon setting $\nu=\epsilon$
\[T_1=\lim_{\epsilon\to 0} \frac{\pi}{2\sin \pi\epsilon} \sum_{m,n\geq0}\frac{(-)^m a^m (\fs b)^{2n}}{m! n!}  \bl\{\frac{\zeta(-2n\!-\!m\!+\!2\epsilon)(\fs b)^{-2\epsilon}}{\g(n\!+\!1\!-\!\epsilon)}-\frac{\zeta(-2n\!-\!m)}{\g(n\!+\!1\!+\!\epsilon)}\br\}\]
\[=\sum_{m,n\geq0}\frac{(-)^m a^m (\fs b)^{2n}}{m! (n!)^2}\bl\{\zeta'(-2n\!-\!m)+\zeta(-2n\!-\!m)\{\psi(n+1)-\log\,\fs b\}\br\}.\]
Upon noting that the hypergeometric function in (\ref{e46}) reduces to an arcsin when $\nu=0$, we then find
\[\sum_{n=1}^\infty e^{-an} K_0(bn)=
\sum_{m,n\geq0}\frac{(-)^m a^m (\fs b)^{2n}}{m! (n!)^2}\bl\{\zeta'(-2n\!-\!m)+\zeta(-2n\!-\!m)\{\psi(n+1)-\log\,\fs b\}\br\}\]
\bee\label{e47}
+\left\{\begin{array}{ll}\dfrac{\fs\pi-\arcsin (a/b)}{b \sqrt{1-(a/b)^2}} & (a<b)\\\dfrac{1}{b} & (a=b). \end{array}\right.
\ee
provided $a<2\pi$ and $b<2\pi$.

Finally, the alternating version of (\ref{e41}) is given by
\[{\hat S}_{K,\nu}(a,b)=\sum_{n=1}^\infty \frac{(-)^{n-1} e^{-an}}{ (\fs bn)^{\nu}}\, K_\nu(bn).\]
It is easily seen that
\[{\hat S}_{K,\nu}(a,b)=S_{K,\nu}(a,b)-S_{K,\nu}(2a,2b),\]
from which the expansion as $a\to 0$ can be obtained from Theorem 3.

\vspace{0.6cm}

\begin{center}
{\bf Appendix: The case of $S_{K,\nu}(a,b)$ when $\nu$ has half-integer values}
\end{center}
\setcounter{section}{1}
\setcounter{equation}{0}
\renewcommand{\theequation}{\Alph{section}.\arabic{equation}}
When $\nu=m+\fs$ we have
\[x^{-m-\frac{1}{2}}K_{m+\frac{1}{2}}(x)=\sqrt{\frac{\pi}{2}}\,\frac{e^{-x}}{x^{m+1}} \sum_{k=0}^m \frac{c_k(m)}{x^k}\qquad (m=0, 1, 2, \ldots),\]
where the $c_k(m)$ are known constants with $c_0(m)=1$. Then $S_{K,\nu}(a,b)$ reduces to
\[S_{K,m+\frac{1}{2}}(a,b)=\frac{2^m\sqrt{\pi}}{b^{m+1}}\,\sum_{k=0}^m \frac{c_k(m)}{b^{k}} \sum_{n=1}^\infty \frac{e^{-(a+b) n}}{n^{m+k+1}}.\]

The sum $\sum_{n\geq 1}n^{-w}e^{-\alpha n}$ for integer $w\geq 1$ has been discussed in \cite[(4.2.10)]{PK}, where it is shown that
\[\sum_{n=1}^\infty n^{-w} e^{-\alpha n}=\frac{(-\alpha)^{w-1}}{\g(w)}\,\{\gamma-\log\,\alpha+\psi(w)\}+\sum_{k=0}^\infty{}^{'} \frac{(-)^k}{k!} \zeta(w-k) \alpha^k\qquad (\alpha<2\pi),\]
where $\gamma$ is the Euler-Mascheroni constant, $\psi(x)$ is the logarithmic derivative of the gamma function and the prime on the summation sign denotes the omission of the term corresponding to $k=w-1$. From this last result it is then possible to construct the expansion of $S_{K,\nu}(a,b)$ when $\nu=m+\fs$. The simplest case when $\nu=\fs$ yields
\bee\label{e48}
S_{K,\frac{1}{2}}(a,b)=\frac{\sqrt{\pi}}{b}\bl\{\sum_{k=1}^\infty \frac{(-)^k}{k!} \zeta(1-k) (a+b)^k-\log\,(a+b)\br\}
\ee
when $a+b<2\pi$.

\end{document}